# A Sparse Differential Algebraic Equation (DAE) and Stiff Ordinary Differential Equation (ODE) Solver in Maple

Taejin Jang[1,=], Maitri Uppaluri[1,=] and Venkat R. Subramanian[1,2,3]

[1]*Materials Science and Engineering Program, Texas Materials Institute, The University of Texas at Austin, Austin, Texas 78712, United States of America*

[2]*Walker Department of Mechanical Engineering, The University of Texas at Austin, Austin, Texas 78712, United States of America*

[3]*Oden Institute for Computational Science and Engineering, The University of Texas at Austin, Austin, TX, United States of America*

[=]Both authors contributed equally to this work.

# Abstract

In this paper, some single-step methods like Trapezoid (TR), Implicit-mid point (IMP), Euler-backward (EB), and Radau IIA (Rad) methods are implemented with adaptive time-stepping in Maple to solve index-1 nonlinear Differential Algebraic Equations (DAEs). Maple's robust and efficient ability to search within a list/set is exploited to identify the sparsity pattern and the analytic Jacobian. The algorithm and implementation were found to be robust and efficient for index-1 DAE problems and scales well for finite difference/finite element discretization of two-dimensional models with system size up to 10,000 nonlinear DAEs and solves the same in a few seconds.

## 1. Introduction and Objectives

The goal of the paper and the code is to solve Hessenberg index-1 DAEs of the form

$$\frac{dy}{dt} = f(y,z)$$
$$0 = g(y,z) \qquad (1)$$

where y is the is a vector of differential variables and z is a vector of algebraic variables. f and g are assumed to be differentiable functions. Note that only autonomous systems are considered in this paper. If the time t explicitly occurs in a system, it can be converted to the form in equation (1) using a dummy variable y defined by dy/dt = 1.

As of today, SUNDIALS implementation of BDF formula (IDA)[1], DASKR[2], and Hairer's FORTRAN implementation of RADAU[3], MATLAB's ode15s[4], Cash's MEBDFDAE[5], Julia's wide range of in-house solvers[6], Maple's inbuilt solvers are good solvers for index 1- DAEs[7].

The goal of this paper is not to arrive at an algorithm that is faster than other implementations in other platforms or codes. The efficiency of a code is a strong function of multiple testing and tuning of adaptive-solver parameters, choice of Jacobian updates, linear solvers, etc. Rather, our goal in this paper is to develop a robust DAE solver in Maple that scales well for large systems. By developing and testing the codes, it turns out that the developed solver and its implementation is competitive compared to many existing solvers. The current implementation can be further optimized by compiling many of the steps involved in C and tuning the solver parameters.

This paper and the code are built on the strengths of Maple, and its ability to do a fast search of variables in a list/set and perform symbolic differentiation and provide for analytic Jacobian. As the system size gets bigger, it is necessary to use sparse linear algebra. The code developed can be used to solve a wide range of index-1 DAEs with a minimal set of inputs from the users and is released under MIT license without any restrictions.

## 2. Algorithms Considered

While the provided code can be easily modified for any single-step method, as a first step, four different algorithms are considered for solving DAEs of the form (1)

**Euler-Backward Method (EB)**
This can be written for equation (1) as

$$y_1 = y_0 + hf(y_1, z_1)$$
$$0 = g(y_1, z_1) \qquad (2)$$

where $y_0$ is the initial condition for y, h is the time step, $y_1$, $z_1$ are ordinary differential equation (ODE) and algebraic variables after the first step, t = h. Note that $z_0$ (initial condition for algebraic variables) does not appear in (2), so any backward-difference method (BDF) method should ideally not need an exact initial condition for z, but code failures can happen if h is too large and if initial guess for z is far off from the expected value when performing a Newton-Raphson iteration to solve equation (2). Typically, f and g are nonlinear functions of y and z. Equation (2) is a system of nonlinear equations of size $N_{ode} + N_{ae}$, where $N_{ode}$ is the number of ODE variables, and $N_{ae}$ is the number of algebraic variables.

**Trapezoid Method (TR or CN)**
This is given by

$$y_1 = y_0 + \frac{h}{2} f(y_1, z_1) + \frac{h}{2} f(y_0, z_0)$$
$$0 = g(y_1, z_1) \tag{3}$$

TR method is also called as Crank Nicolson (CN) method for solving semi-discretized systems of parabolic partial differential equations (PDEs) and will be referred to as CN in this paper and code. For CN method, there is a need to solve $g(y_0,z_0)$ to find a consistent initial condition at t = 0. The code is written such that $z_0$ is found for all the algorithms before marching in time, as small errors in z can affect the results in future times. This is particularly important for state estimation as needed for Battery-Management Systems (BMS) in battery modeling and simulation.

**Implicit Mid-Point-Trapezoid Method (IMPTRAP)**
The implicit-midpoint algorithm is given by

$$y_1 = y_0 + hf(y_{mid}, z_{mid})$$
$$0 = g(y_{mid}, z_{mid}) \tag{4}$$

In the implicit-mid point method, residues are evaluated at the mid-point of the time interval h/2. Since the method is not stiffly accurate, there is a need to find $z_1$ by solving $g(y_1,z_1)$ after finding $y_{mid}$ and $z_{mid}$. $y_{mid}$ can be approximated using linear-interpolation as $y_{mid} = (y_1+y_0)/2$. This does not introduce additional errors.

$$y_1 = y_0 + hf\left(\frac{(y_0 + y_1)}{2}, z_{mid}\right)$$
$$0 = g\left(\frac{(y_0 + y_1)}{2}, z_{mid}\right) \tag{5}$$

Equation (5) still requires a nonlinear solution for $z_1$ at the end of each step. In this paper, $z_{mid}$ is approximated to be $z_{mid} = (z_0+z_1)/2$ for the efficiency of implementation and $g(y_1,z_1)$ is forced to be zero (as a projected step). This method can be viewed as a hybridization of the implicit-midpoint method for ODE variables and the trapezoid method for algebraic variables. This hybridization approach seems to work well for smooth DAEs, and it might fail for problems in which z is discontinuous with time or for very stiff DAEs that require L-stable schemes.[8]

$$y_1 = y_0 + hf\left(\frac{(y_0+y_1)}{2}, \frac{(z_0+z_1)}{2}\right) \tag{6}$$

$$0 = g(y_1, z_1)$$

**RadauIIA Method (Rad)**

The third-order RadauIIA method for equation (1) can be written as

$$\begin{aligned} y_1 &= y_0 + \frac{h}{4}f(y_1,z_1) + \frac{3h}{4}f(y_{int},z_{int}) \\ 0 &= g(y_1,z_1) \\ y_{int} &= y_0 - \frac{h}{12}f(y_1,z_1) + \frac{5h}{12}f(y_{int},z_{int}) \\ 0 &= g(y_{int},z_{int}) \end{aligned} \tag{7}$$

Where $y_{int}$ and $z_{int}$ are ODE and algebraic variables at t = h/3. The system size is doubled for the Newton-Raphson solver for this method.

Equation (7) can be rewritten (by using the inverse of Runge-Kutta coefficient matrix A) as

$$\begin{aligned} \frac{5}{2}(y_1-y_0) - \frac{9}{2}(y_{int}-y_0) &= hf(y_1,z_1) \\ 0 &= g(y_1,z_1) \\ \frac{1}{2}(y_1-y_0) + \frac{3}{2}(y_{int}-y_0) &= hf(y_{int},z_{int}) \\ 0 &= g(y_{int},z_{int}) \end{aligned} \tag{8}$$

This change makes a significant difference in computational efficiency. For finite difference discretization of a single PDE in 2D with a five-point stencil, equation (7) will have 10 non-zero entries in the Jacobian for a particular row as opposed to only 6 non-zero entries in equation (8).

All the methods require a nonlinear solver to update and find $y_1$ and $z_1$. To facilitate this, difference between y and z are defined as uu

$$y_1 = uu_i + y_0, \; i = 1..N_{ode}, \; z_1 = uu_i + z_0, \; i = N_{ode}+1..N_{ode}+N_{ae}. \tag{9}$$

$y_1$ and $z_1$ are stacked together in vector form as $\mathbf{Y_1}$

$$\mathbf{Y_1} = \mathbf{Y_0} + \mathbf{uu} \tag{10}$$

The four algorithms are written in residual form (after defining $N_t = N_{ode}+N_{ae}$) as

$$uu_i - h \cdot f_i(\mathbf{uu}+\mathbf{Y_0}) = 0, i = 1..N_{ode}$$
$$g_i(\mathbf{uu}+\mathbf{Y_0}) = 0, i = N_{ode}+1.. N_{ode}+N_{ae} \quad (11)\ (EB)$$

$$uu_i - \frac{h}{2} f_i(\mathbf{uu}+\mathbf{Y_0}) - \frac{h}{2} f_i(\mathbf{Y_0}) = 0, i = 1..N_{ode}$$
$$g_i(\mathbf{uu}+\mathbf{Y_0}) = 0, i = N_{ode}+1.. N_{ode}+N_{ae} \quad (12)\ (CN)$$

$$uu_i - h \cdot f_i\left(\frac{\mathbf{uu}}{2}+\mathbf{Y_0}\right) = 0, i = 1..N_{ode}$$
$$g_i(\mathbf{uu}+\mathbf{Y_0}) = 0, i = N_{ode}+1.. N_{ode}+N_{ae} \quad (13)\ (IMP)$$

$$\frac{5}{2}uu_i - \frac{9}{2}uu_{i+N_{ode}+N_{ae}} - h \cdot f_i(\mathbf{uu}+\mathbf{Y_0}) = 0, i = 1..N_{ode}$$
$$g_i(\mathbf{uu}+\mathbf{Y_0}) = 0, i = N_{ode}+1.. (N_{ode}+N_{ae})$$
$$\frac{1}{2}uu_{i-(N_{ode}+N_{ae})} + \frac{3}{2}uu_i - h \cdot f_i(\mathbf{uu}+\mathbf{Y_0}) = 0, i = (N_{ode}+N_{ae})+1..(N_{ode}+N_{ae})+N_{ode}$$
$$g_i(\mathbf{uu}+\mathbf{Y_0}) = 0, i = (N_{ode}+N_{ae})+N_{ode}+1 .. 2(N_{ode}+N_{ae}) \quad (14)\ (Rad)$$

All the methods involve finding the **uu** vector for a given **Y₀** vector as input. All the methods can be implemented with h = 0 to find the consistent initial condition for z at t = 0 without creating a new set of equations just for initialization.

## 3. Error-Control and Time-Stepping

Error control in the code is achieved based on the absolute tolerance requirement. In this paper, relative tolerance is taken to be a scalar quantity and set to 10 times the absolute tolerance. This is achieved by running a single step of the algorithm with time step= h, and two-half steps (t=h/2 twice to complete the same step) and finding the error using two different estimates.

$$y_{err} = \frac{(y_{h_2} - y_h)}{(2^p - 1)}$$
$$err = \left\| \frac{y_{err}[i]}{(a_{tol} + y_{err}[i] \cdot r_{tol})} \right\| \quad (15)$$
$$h_{new} = \min\left(h_{max}, h_{old} \min\left(3.0, 0.9\left(\frac{1}{err}\right)^{\frac{1}{(p+1)}}\right)\right)$$

In equation (15), p is the order of accuracy, which is 1,2,2,3 for EB, CN, IMPTRAP and Rad methods considered in this paper. Both the ODE and algebraic variables are included in the error-estimates. While many single-step methods can have cheaper error estimates using embedded methods, our experience suggests that this half-step approach provides a more robust error estimate for the algebraic variables. In addition, the approach also provides a higher order accuracy at the end of each time step with Richardson extrapolation as

$$y_R = \frac{\left(2^p y_{h_2} - y_h\right)}{2^p - 1} \tag{16}$$

For certain stiff problems, this might cause stability issues and the code can be modified to use $y_{h2}$ at the end of each step. Also, both ODE and algebraic variables are extrapolated using this formula at the end of each step.

While the number of Newton Raphson iterations required for convergence can be used to find a criterion for updating the Jacobian, in this code if err is greater than 0.1, the Jacobian is updated. If err is greater than 1, the step is rejected and the time step is reduced by 4. Also, the Jacobian obtained at a particular time t is used for all the linear solves for the Newton Raphson method for every iteration for both the calculation with t = h and two repeated steps with t = h/2. This means that the code will not perform LU Decomposition more than once for a given time step.

## 4. Symbolic Math and Use of Search Tools

One of the unique aspects of the developed code is the analytic Jacobian found using search tools in Maple. When a wide range of spatial discretization methods is used to convert PDEs to DAEs (finite difference, collocation, spectral, finite volume methods, *etc.*), the sparsity pattern is hard to identify and code for general PDEs and boundary conditions for different spatial discretization approaches. In this code, for a given system of DAEs, Maple's search and symbolic capabilities are exploited to arrive at a robust and efficient way to (1) search, sort and label variables and indices (2) differentiate expressions for analytic Jacobian (3) create sparse Jacobian and procedures for the same.

The code implemented scales well for two-dimensional PDEs and takes very little time (seconds) to find the analytic sparse Jacobian even for >100,000 DAEs resulting from semi-discretization of two-dimensional PDEs.

The specific algorithm in Maple that performs this search and stores the sparse Jacobian is given below.

```
for i to Nt do
    L := indets(Equation[i]);
    LL := [seq(ListTools:-Search(L[j], Vars1), j = 1 .. nops(L))];
    LL2 := ListTools:-MakeUnique(LL);
    if LL2[nops(LL2)] = 0 then
        LL := [seq(LL2[i], i = 1 .. nops(LL2) - 1)];
    end if;
    if LL2[1]= 0 then LL:=[seq(LL2[i],i=2..nops(LL2))]: end:
    for j to nops(LL) do
```

```
            j11[i, LL[j]] := diff(Equation[i], uu[LL[j]]);
          end do;
        end do;
```
where Nt is the total number of variables, all the residual equations are stored in Equation, j11 is the analytic sparse Jacobian.

## 5. Using the Codes

The solver can be easily called from Maple to solve index-1 DAEs. This is best explained by demonstrating code use for a simple example.

**Example 1:**
Consider the DAE system

$$\frac{dy}{dt} = z$$
$$y^2 + x^2 = 1 \qquad (17)$$
$$y(0) = 1; z(0) = 0.95$$

The exact initial condition for z is 1.0, but 0.95 is provided as an approximate value for z to test the code's ability to initialize algebraic variables. First, the solver is called from Maple, and then the equations are called and solved as below.

```
> read("DAESolver.txt"):
> eq1:=diff(y(t),t)=z(t);eq2:=y(t)^2+z(t)^2-1=0;
```
$$eq1 := \frac{d}{dt} y(t) = z(t)$$
$$eq2 := y(t)^2 + z(t)^2 - 1 = 0$$
```
> ICs:=[y(t)=0.,z(t)=0.95];
```
$$ICs := [y(t) = 0., z(t) = 0.95]$$
```
> eqodes:=[eq1];eqaes:=[eq2];
```
$$eqodes := \left[\frac{d}{dt} y(t) = z(t)\right]$$
$$eqaes := [y(t)^2 + z(t)^2 - 1 = 0]$$
```
> solveroptions:=[tf=1.0,atol=1e-4,hinit=1e-
5,hmax=0.05,Ntot=1000,iter=5];#Don't change the order
```
$$solveroptions := [tf = 1.0, atol = 0.0001, hinit = 0.00001,$$
$$hmax = 0.05, Ntot = 1000, iter = 5]$$
```
> sol1:=IMPDAE(eqodes,eqaes,ICs,solveroptions);
```

$$sol1 := [[[t = 0., y(t) = 0., z(t) = 1.00000000000000]],$$
$$[[t = 0.00001, y(t) = 0.00000999999999983330, z(t)$$
$$= 0.999999999949997]],$$
$$[[t = 0.000040, y(t) = 0.0000399999999893332, z(t)$$
$$= 0.999999999199997]],$$
$$[[t = 0.0001300, y(t) = 0.000129999999633833, z(t)$$
$$= 0.999999991549997]],$$
$$[[t = 0.00040000, y(t) = 0.000399999989333331,$$
$$z(t) = 0.999999919999998]],$$
$$[[t = 0.001210000, y(t) = 0.00120999970473985,$$
$$z(t) = 0.999999267950086]],$$
$$[[t = 0.0036400000, y(t) = 0.00363999196191437,$$
$$z(t) = 0.999993375207312]],$$
$$[[t = 0.01093000000, y(t) = 0.0109297823760090,$$
$$z(t) = 0.999940268144656]],$$
$$[[t = 0.032800000000, y(t) = 0.0327941190418480,$$
$$z(t) = 0.999462128225188]],$$
$$[[t = 0.082800000000, y(t) = 0.0827054207805036,$$
$$z(t) = 0.996574038074572]],$$
$$["29\ element\ Vector[column]"]]$$

**> ntot:=Dimension(sol1);**

$$ntot := 29$$

**>
plot([seq(subs(sol1[i],[t,z(t)]),i=1..ntot)],style=point,thickness=3,axes=boxed,labels=[t,y(t)]);**

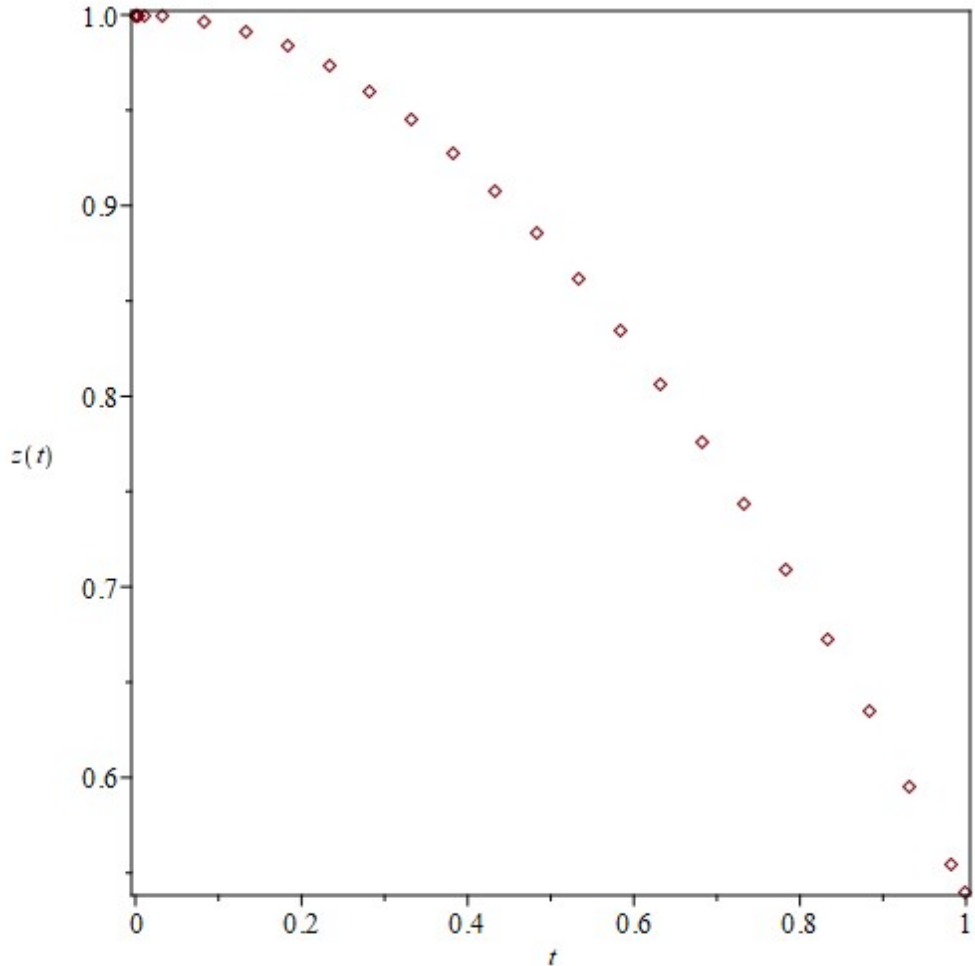

*Figure 1: z vs t plot using the data obtained from the simulation of example 1 with the implemented DAE solver with the IMP algorithm. The code is able to initialize z correctly and then solve the DAE model.*

The solver took only 29 time-steps to simulate this example. The user has to stack all the ODEs in eqodes and all the AEs in eqaes. If there are no algebraic variables, the users should provide an empty list [] for eqaes.

The solver options include the following inputs

$t_f$ = total time of simulation

$a_{tol}$ = absolute tolerance expected for the simulation

$h_{init}$ = starting time step for simulation, expected to be less than tf. Suggested value is $h_{init}$ = $\min(10^{-6}, t_f.a_{tol})$

$h_{max}$ = maximum value of the time step. If $h_{max}$ is too large, plots may not include enough data points and for some problems, code might fail if $h_{max}$ is too large. Suggested value in the code is $h_{max} = t_f/20$,

$N_{tot}$ = maximum number of time steps, taken to be 1000.

The initial conditions for ODE and algebraic variables are listed in ICs. This should follow the order specified in eqodes (list of ODEs) and then eqaes (list of AEs). The initial conditions are provided for algebraic variables as well. These can be guesses for the algebraic variables as the solver does the initialization. While Maple can automatically find and separate ODEs and AEs from a given system, this is avoided in the solver for efficiency purposes.

The solver can be called to provide the CPU time as well.

**Example 2:**

Next, Vander Pol's system is considered. This is a test case for stiff solvers, and the developed solver is able to simulate this model as shown in Figure 2.

$$\frac{dx}{dt} = \mu(1-y^2)x - y$$
$$\frac{dy}{dt} = x \qquad (18)$$
$$x(0) = 0;\ y(0) = 2;\ \mu = 2$$

```
> eqodes :=[diff(x(t),t)=mu*(1-y(t)^2)*x(t)-
y(t),diff(y(t),t)=x(t)] ;Ics :=[x(t)=0.0,y(t)=2.0] ;mu :=2 ;
```

$$eqodes := \left[\frac{d}{dt}x(t) = \mu\left(1 - y(t)^2\right)x(t) - y(t),\ \frac{d}{dt}y(t) = x(t)\right]$$

$$ICs := [x(t) = 0.,\ y(t) = 2.0]$$

$$\mu := 2$$

```
> soln:=dsolve({op(eqodes),x(0)=0,y(0)=2.0},type=numeric,maxf
un=0,abserr=1e-10):soln(10.0);plots:-
odeplot(soln,[t,y(t)],0..10):
```

$$[t = 10.0,\ x(t) = -1.08904705944854,\ y(t) = 0.841554375093746]$$

$$\mu := 2$$

```
> solveroptions:=[tf=10.0,atol=1e-5,hinit=1e-
6,hmax=0.1,Ntot=1000,iter=5];#Don't't'change the order
```

$$solveroptions := [tf = 10.0,\ atol = 0.00001,\ hinit = 0.000001,\ hmax = 0.1,\ Ntot = 1000,\ iter = 5]$$

```
>
sol21:= CodeTools:- Usage(IMPDAE(eqodes,[],ICs,solveroptions));n
tot:= Dimensions(sol21);s ol21[ntot];s oln(10.0);
```

"Integration completed, but number of failed steps=", 19

```
memory used=68.38MiB, alloc change=-4.00MiB, cpu time=719.00ms, real
time=698.00ms, gc time=187.50ms
```

$sol21 := [[[t=0., x(t) = 0., y(t) = 2.]],$
$[[t=0.000001, x(t) = -0.00000199999400001166,$
$y(t) = 1.99999999999899]],$
$[[t=0.0000040, x(t) = -0.00000799990400074661,$
$y(t) = 1.99999999998399]],$
$[[t=0.00001300, x(t) = -0.0000259989860256310,$
$y(t) = 1.99999999983098]],$
$[[t=0.000040000, x(t) =$
$-0.0000799904007466324, y(t)$
$= 1.99999999840010]],$
$[[t=0.0001210000, x(t) =$
$-0.000241912174665421, y(t) = 1.99999998536251]$
$],$
$[[t=0.00036400000, x(t) =$
$-0.000727205586438103, y(t) = 1.99999986760037]$
$],$
$[[t=0.001093000000, x(t) =$
$-0.00217884732120933, y(t) = 1.99999880795830]$
$],$
$[[t=0.0032800000000, x(t) =$
$-0.00649585978274598, y(t) = 1.99998931183839]$
$],$
$[[t=0.00984100000000, x(t) =$
$-0.0191119252124828, y(t) = 1.99990503369744]],$
$["172 \text{ element Vector[column]"}]]$

$ntot := 172$

$[t=10.0000000000000, x(t) = -1.08901427886770, y(t)$
$= 0.841578964658030]$

$[t=10.0, x(t) = -1.08904705944854, y(t)$
$= 0.841554375093746]$

> 
```
plot([seq(subs(sol21[i],[t,y(t)]),i=1..ntot)],style=point,thickness=3,axes=boxed,labels=[t,y(t)]);
```

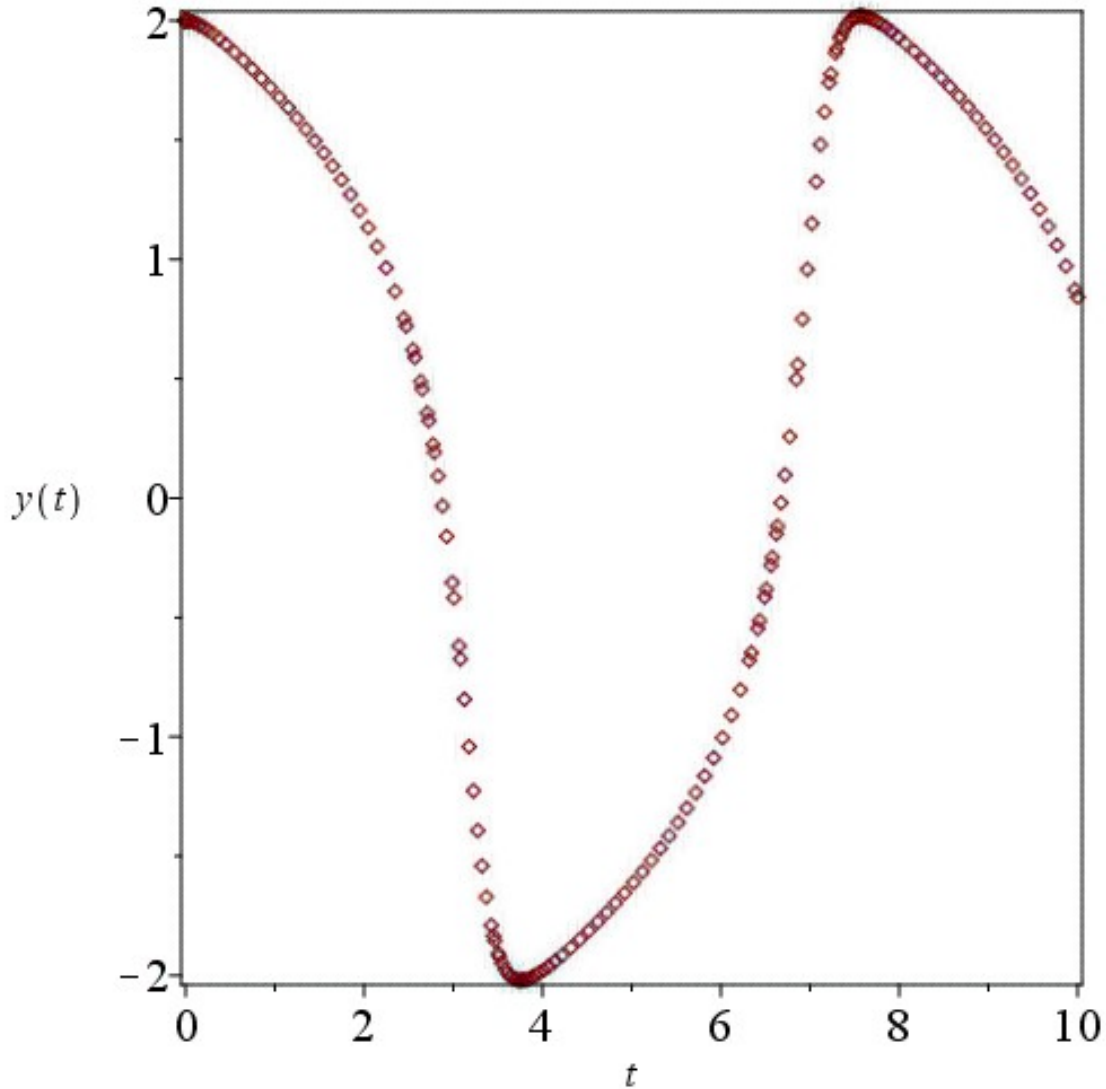

*Figure 2: z vs t plot using the data obtained from the simulation of example 2 (Vander Pol's equation) with the implemented DAE solver and the IMP algorithm. The code is able to solve the stiff system of ODEs.*

**Example 3:**
Next, a simple DAE system which fails with Maple's default DAE solver (because of unknown initial condition for the algebraic variable) is solved. The solver can find the exact initial condition for algebraic variable z and solve the system as shown in Figures 3 and 4.

$$\frac{dy}{dt} = -2y + z^2$$
$$-100\ln(z) + 2y = 5 \tag{19}$$
$$y(0) = 2; z(0) = 1$$

```
> eq1:= diff(y(t),t)=-2*y(t)+z(t)^2;e q2:= -
100*log(z(t))+2*y(t)=5;
```

$$eq1 := \frac{d}{dt} y(t) = -2 y(t) + z(t)^2$$

$$eq2 := -100 \ln(z(t)) + 2 y(t) = 5$$

> **eqodes:=[eq1];eqaes:=[eq2];**

$$eqodes := \left[ \frac{d}{dt} y(t) = -2 y(t) + z(t)^2 \right]$$

$$eqaes := [-100 \ln(z(t)) + 2 y(t) = 5]$$

> **ICs:=[y(t)=2.,z(t)=1.0];**

$$ICs := [y(t) = 2., z(t) = 1.0]$$

> **soln:=dsolve({eq1,eq2,op(subs(t=0,ICs))},type=numeric):**

Error, (in dsolve/numeric/DAE/checkconstraints) the initial conditions do not satisfy the algebraic constraints
  error = 1., tolerance = .100e-3, constraint = -100*ln(z(t))+2*y(t)-5

> **solveroptions:=[tf=5.0,atol=1e-5,hinit=1e-6,hmax=0.5,Ntot=1000,iter=5];#Don't't'change the order**

$$solveroptions := [tf = 5.0, atol = 0.00001, hinit = 0.000001, hmax = 0.5, Ntot = 1000, iter = 5]$$

>
**sol1:= CodeTools:- Usage(IMPDAE(eqodes,eqaes,ICs,solveroptions));**

memory used=11.68MiB, alloc change=0 bytes, cpu time=94.00ms, real time=89.00ms, gc time=0ns

$$sol1 := [[[t = 0., y(t) = 2., z(t) = 0.990049833749168]$$
$$],$$
$$[[t = 0.000001, y(t) = 1.99999698020163, z(t) = 0.990049773954149]],$$
$$[[t = 0.0000040, y(t) = 1.99998792084205, z(t) = 0.990049594569828]],$$
$$[[t = 0.00001300, y(t) = 1.99996074308307, z(t) = 0.990049056423389]],$$
$$[[t = 0.000040000, y(t) = 1.99987921268375, z(t) = 0.990047442042806]],$$
$$[[t = 0.0001210000, y(t) = 1.99963464738218, z(t) = 0.990042599429625]],$$
$$[[t = 0.00036400000, y(t) = 1.99890118449167, z(t) = 0.990028076346013]],$$
$$[[t = 0.001093000000, y(t) = 1.99670289150694, z(t) = 0.989984549867361]],$$
$$[[t = 0.0032800000000, y(t) = 1.99012683489093, z(t) = 0.989854354540289]],$$
$$[[t = 0.00984100000000, y(t) = 1.97056702229335, z(t) = 0.989467202957876]],$$
$$[\text{"33 element Vector[column]"}]]$$

> **ntot:=Dimensions(sol1);**

$$ntot := 33$$

```
> plot([seq(subs(sol1[i],[t,y(t)]),i=1..ntot)],style=point,thickne
ss=3,axes=boxed,labels=[t,y(t)]);
```

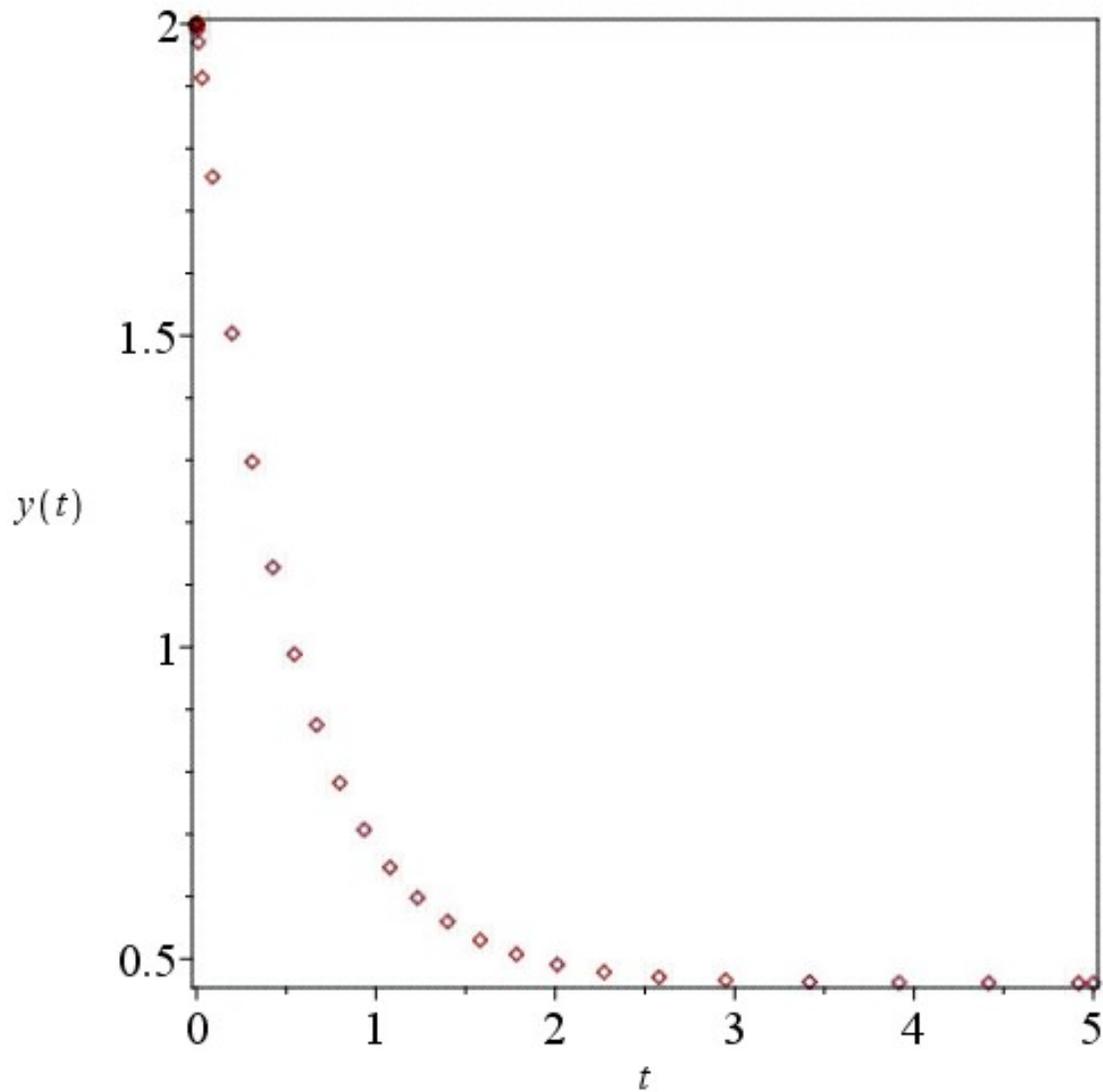

*Figure 3: y (ODE variable) vs t plot using the data obtained from the simulation of example 3 with the implemented DAE solver and the IMP algorithm. The code is able to initialize z correctly and then solve the DAE model.*

```
> plot([seq(subs(sol1[i],[t,z(t)]),i=1..ntot)],style=point,thickne
ss=3,axes=boxed,labels=[t,z(t)]);
```

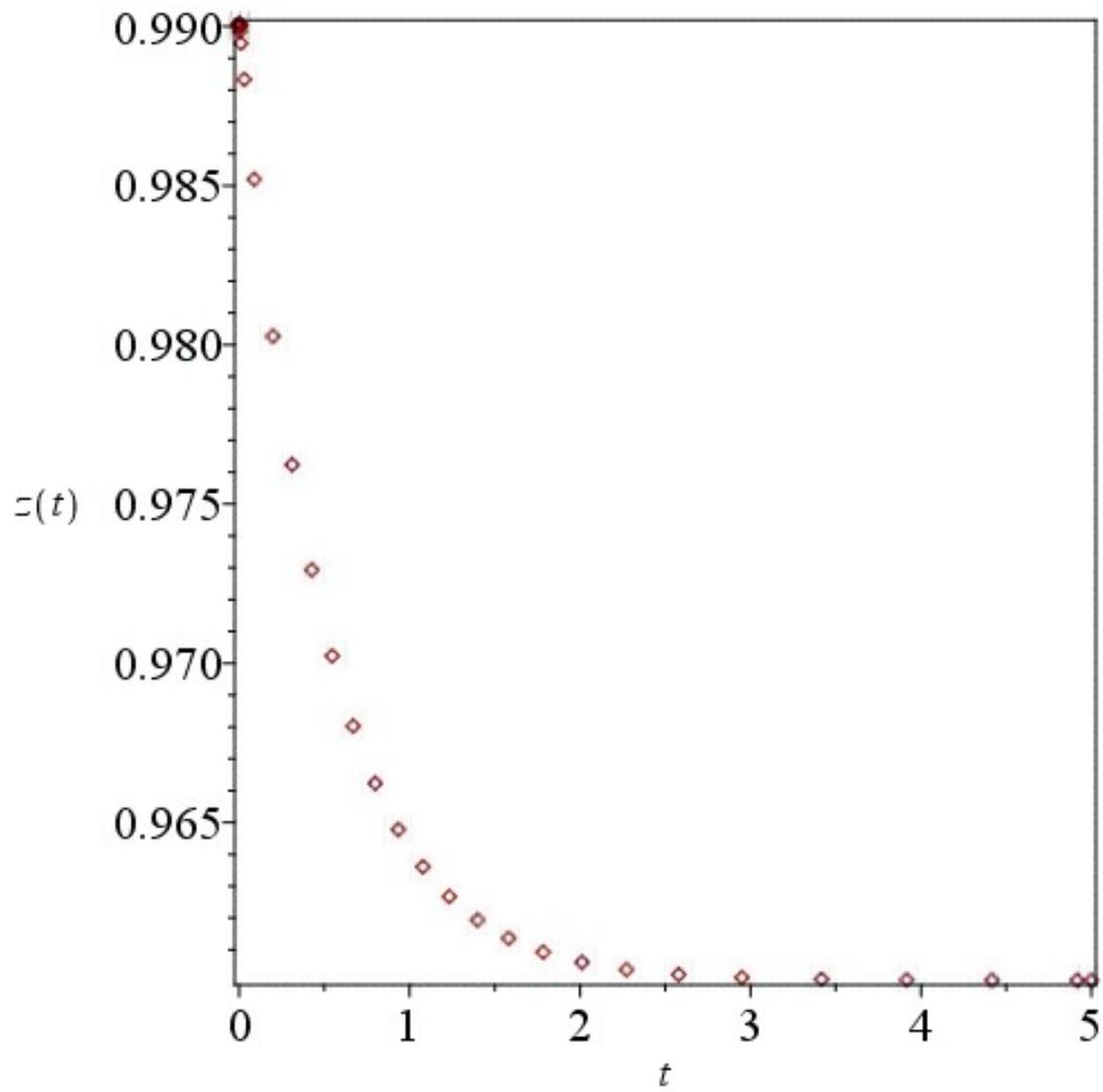

*Figure 4: z (Algebraic variable) vs t plot using the data obtained from the simulation of example 3 with the implemented DAE solver and the IMP algorithm. The code is able to initialize z correctly and then solve the DAE model.*

**Example 4:**
Next, a system of partial differential equations in one-dimension is considered

$$\frac{\partial c}{\partial t} = \frac{\partial^2 c}{\partial x^2} - c(1+z)$$

$$0 = \frac{\partial^2 z}{\partial x^2} - (1-c^2)\exp(-z)$$

$$c(x,0) = 1; z(x,0) = 0 \quad (20)$$

$$\frac{\partial c}{\partial x}(0,t) = 0; \frac{\partial z}{\partial x}(0,t) = 0$$

$$c(1,t) = 1; z(1,t) = 0$$

Discretized form of equation (20) with the cell-centered finite difference approach in x can be written as

$$\frac{dc_i}{dt} = \frac{c_{i+1} - 2c_i + c_{i-1}}{(\Delta x)^2} - c_i(1+z_i), i = 1..N$$

$$0 = \frac{z_{i+1} - 2z_i + z_{i-1}}{(\Delta x)^2} - (1-c_i^2)\exp(-z_i), i = 1..N$$

$$c_i(0) = 1; z_i(0) = 0, i = 1..N \quad (21)$$

$$\frac{c_1 - c_0}{\Delta x} = 0; \frac{z_1 - z_0}{\Delta x} = 0$$

$$\frac{c_N + c_{N+1}}{2} = 1; \frac{z_N + z_{N+1}}{2} = 0$$

where N is the number of elements in x.

The code is written in such that N can be changed from 2,4,8,16 until satisfactory results are seen. It is noted that the code works even for N = 10,000 node points in x in very few seconds. For this problem, c and z at x = 0 and t = 1 (given by $(c_0+c_1)/2$ and $(z_0+z_1)/2$) converge with increasing values of N 4,8,16,32,64 as
[0.708501773693253, -0.276198079090988],
[0.706802805832433, -0.274160215413445],
[0.706377861591221, -0.273646880455202],
[0.706273290106072, -0.273518340174281],
[0.706246833067691, -0.273486190782169].

The CPU time for all the examples is tabulated in Table 2 for all the methods using Intel® Core™ i9–12900K CPU and 32 GB RAM.

**Example 5:**
A 2D problem is next considered

$$\frac{\partial c}{\partial t} = \frac{\partial^2 c}{\partial x^2} + \frac{\partial^2 c}{\partial y^2} - \phi^2 c^2$$
$$c(x, y, 0) = 0$$
$$\frac{\partial c}{\partial x}(0, y, t) = 0; \frac{\partial c}{\partial y}(x, 0, t) = 0 \quad (22)$$
$$c(x, 1, t) = 1; c(1, y, t) = 1$$

$\phi$ can be varied as 0.1, 1, and 10 to see the need for different number of node points in x. A cell-centered finite difference (FD) was used to simulate this model for $\phi = 0.5$. Since the Maple code is provided and FD scheme for a more detailed model is provided for example 6, the FD scheme is not provided for this model.

**Example 6:**
Next, the concentration and potential distribution (tertiary current-distribution) in an electrolyte (electrochemical system) are considered. With the dilute solution theory and the electroneutrality assumption, the model considered is

$$\frac{\partial c}{\partial t} = D_x \frac{\partial^2 c}{\partial x^2} + D_y \frac{\partial^2 c}{\partial y^2}$$
$$\frac{\partial}{\partial x}\left(D_x c \frac{\partial \phi}{\partial x}\right) + \frac{\partial}{\partial y}\left(D_y c \frac{\partial \phi}{\partial y}\right) = 0 \quad (23)$$

with the boundary conditions

$$\left.\frac{\partial c}{\partial y}\right|_{y=0} = 0, \left.\frac{\partial \phi}{\partial y}\right|_{y=0} = 0$$

$$\left.\frac{\partial c}{\partial y}\right|_{y=H} = 0, \left.\frac{\partial \phi}{\partial y}\right|_{y=H} = 0$$

$$D_x \left.\frac{\partial c}{\partial x}\right|_{x=L} = \delta, D_x c \left.\frac{\partial \phi}{\partial x}\right|_{x=L} = \delta \quad (24)$$

$$\begin{cases} D_x \left.\frac{\partial c}{\partial x}\right|_{x=0} = Da \cdot c \cdot \phi, \ D_x \left.\frac{\partial \phi}{\partial x}\right|_{x=0} = Da \cdot \phi & 0 < y <= \frac{H}{2} \\ \left.\frac{\partial c}{\partial x}\right|_{x=0} = 0, \ \left.\frac{\partial \phi}{\partial x}\right|_{x=0} = 0 & \frac{H}{2} < y <= H \end{cases}$$

The initial condition for c was taken to be 1 everywhere. Equation (23) is the final form of the model for the transport of a binary electrolyte based on Nernst-Planck equation for diffusive and migrative flux, coupled with electroneutrality. Equal diffusivities were assumed for both the cation and anion and the model assumes different constant diffusivities in the x and y direction.[8]

Discretized form of equation (23) with cell-centered finite difference method in x and y can be written as

$$\frac{dc_{i,j}}{dt} = D_x \frac{c_{i+1,j} - 2c_{i,j} + c_{i-1,j}}{(\Delta x)^2} + D_y \frac{c_{i,j+1} - 2c_{i,j} + c_{i,j-1}}{(\Delta y)^2}, i=1..N, j=1..M$$

$$0 =$$

$$\frac{\left(D_x \frac{c_{i+1,j} + c_{i,j}}{2} \left(\frac{\phi_{i+1,j} - \phi_{i,j}}{\Delta x}\right)\right) - \left(D_x \frac{c_{i,j} + c_{i-1,j}}{2} \left(\frac{\phi_{i,j} - \phi_{i-1,j}}{\Delta x}\right)\right)}{\Delta x}$$

$$+ \frac{\left(D_y \frac{c_{i,j} + c_{i,j+1}}{2} \left(\frac{\phi_{i,j+1} - \phi_{i,j}}{\Delta y}\right)\right) - \left(D_y \frac{c_{i,j} + c_{i,j-1}}{2} \left(\frac{\phi_{i,j} - \phi_{i,j-1}}{\Delta y}\right)\right)}{\Delta y}, i=1..N, j=1..M$$

$$\frac{c_{i,1} - c_{i,0}}{\Delta y} = 0; \frac{\phi_{i,1} - \phi_{i,0}}{\Delta y} = 0, i = 1..N$$

$$\frac{c_{i,M+1} - c_{i,M}}{\Delta y} = 0; \frac{\phi_{i,M+1} - \phi_{i,M}}{\Delta y} = 0, i = 1..N$$

$$D_x \frac{c_{N+1,j} - c_{N+1,j}}{\Delta x} = \delta; D_x \frac{c_{N+1,j} + c_{N,j}}{2} \left( \frac{\phi_{N+1,j} - \phi_{N,j}}{\Delta x} \right) = \delta, j = 1..M \qquad (25)$$

$$\begin{cases} D_x \frac{c_{1,j} - c_{0,j}}{\Delta x} = Da \left( \frac{c_{0,j} + c_{1,j}}{2} \right) \left( \frac{\phi_{0,j} + \phi_{1,j}}{2} \right); D_x \frac{\phi_{1,j} - \phi_{0,j}}{\Delta x} = Da \left( \frac{c_{0,j} + c_{1,j}}{2} \right) = 0, j = 1..M/2 \\ \frac{c_{1,j} - c_{0,j}}{\Delta x} = 0; \frac{\phi_{1,j} - \phi_{0,j}}{\Delta x} = 0, j = M/2 + 1..M \end{cases}$$

In equation (25), Da is The Damkohler number, δ is the applied current density in dimensionless form. Typically, either Dx or Dy can be eliminated (with scaling), but we leave it as is in equation 25 and the code.

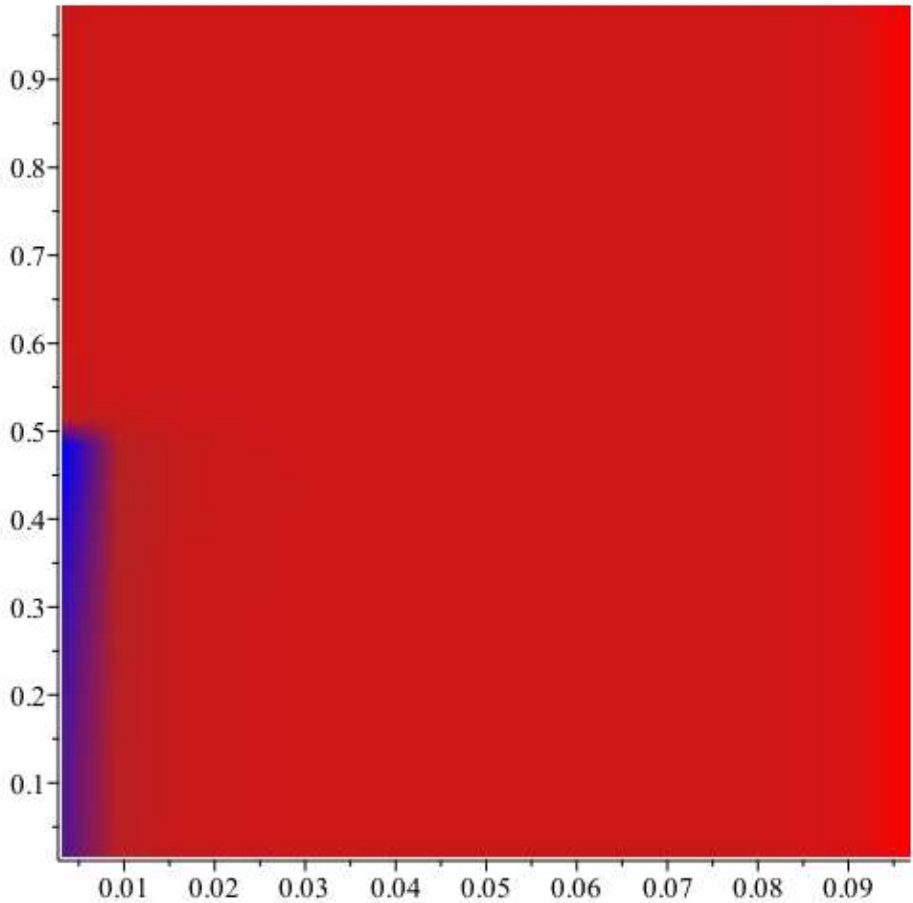

Figure 5: *Surface plot of concentration at short times for the electrochemical model. A boundary layer is formed near the electrode at x =0 which then diffuses with time*

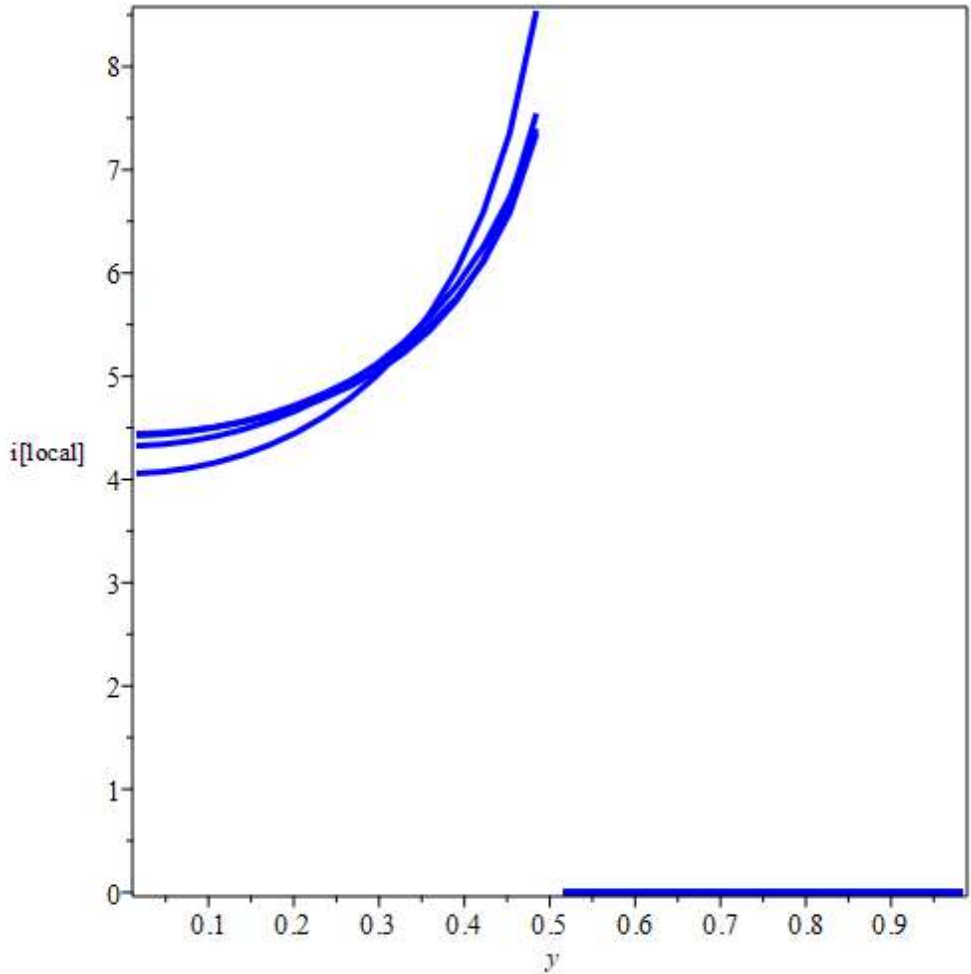

Figure 6: *Current density (flux distribution) at x = 0 from the tertiary current distribution model at different times. Singularity is seen at y = 0.5 which reduces the observed order spatial discretization accuracy.*

Results at x =0.5, y=0 and x=0, y = 0.5 are used to study the convergence of results for N with M = 2N, $\Delta x$ = 0.1/N, $\Delta y$ = 1/M. Results obtained for different values of N= [4,8,16,32,64,128,256] are

[0.731415952972219, 0.901537165192901, 0.827803387614459, 1.003133599893680],
[0.738562409389896, 0.868922997011145, 0.817870608768577, 0.962556079339768],
[0.741015294469612, 0.858571490970910, 0.813050089056554, 0.948343507838869],
[0.741814512872462, 0.855333932557370, 0.811016709652655, 0.943505084693148],
[0.742062660273571, 0.854350929972723, 0.810241580007608, 0.941922913444259],
[0.742136965090761, 0.854060601152856, 0.809966832237226, 0.941423593071616].

One can see that the code scales well for large number of node points which are needed for this problem because of singularity at y = 0.5 and x = 0. Use of variable grids in x and y and other numerical methods can help in solving the problem more efficiently. Convergence analysis and finding the optimal method for spatial discretization are beyond the scope of this paper and the code developed. Surface plots at short times and current distribution at x = 0 are also plotted in figures 5 and 6.

## 6. Analysis of Algorithms and Results

A comparison of properties of different algorithms considered in this paper is reported below in Table 1. Computational results are summarized in Table 2. For low tolerances (high accuracies), RadauIIA method is recommended for a small number of DAEs. For a high number of DAEs, linear solver dominates the total computational cost. For the PDE examples discussed in this paper, IMPTRAP method was found to be the most efficient for an absolute tolerance of $10^{-6}$. Both Euler-backward and RadauIIA methods are L-stable and are recommended for very stiff DAEs. For Hamiltonian DAEs, IMPTRAP method offers symplecticity for the ODE variables only as the Trapezoid method is only weakly symplectic (symmetric). For problems requiring monotonicity, only Euler-backward method is offered.

Table 1. Properties of different time-stepping algorithms implemented in this paper.

| Method | Order | A & L-Stability | Extrapolated Order |
|---|---|---|---|
| Euler Backward | 1 | Yes & Yes | 2 |
| Crank-Nicolson | 2 | Yes & No | 4 |
| Implicit mid-point-Trapezoid | 2 | Yes & No | 4 |
| RadauIIA | 3 | Yes & Yes | 4 |

Additional comments about different examples are given below.

Example 1: This is a simple index-1 DAE, but the profile is oscillatory in nature. The code presented in this paper is able to meet the tolerances specified. For a set absolute tolerance of $10^{-6}$, the proposed solver is able to simulate this model in 33 time-steps and CPU time of 58 ms with the Radau method with 2 failed time steps. Note that max time step was chosen to be 0.05, otherwise the code might run even faster.

Example 2: Vander Pol's model is a stiff system of ODEs. So, an empty list is passed to the solver for the algebraic equations. For a set absolute tolerance of $10^{-6}$, the solver is able to simulate this model in 114 time-steps and CPU time of 293 ms with the Radau method with 17 failed time steps. Note that max time step was chosen to be 0.1, otherwise the code might run even faster. One can see that the solver is able to adapt and use smaller time-steps as needed.

Example 3: This simple DAE system was chosen to show the importance of consistent initial conditions for the algebraic equations and variables. Maple's inbuilt DAE solver fails to solve this DAE for the initial condition of z =1 for y =2. The implemented solver is able to simulate this model in 34 time-steps and CPU time of 55 ms with the Radau method with 0 failed time steps. Note that max time step was chosen to be 0.5, otherwise the code might run even faster.

Example 4: This example was chosen to show the scalability of the solver for a large-scale system and for highly nonlinear algebraic equations. The PDE system is simulated with a cell-centered finite difference method in x. Maple's inbuilt DAE solver fails to solve PDE system for N>5 elements in x. The implemented solver is able to simulate this model in 45 time-steps and CPU time of 311 ms with the IMPTRAP method with 0 failed time-steps with N= 128 elements (total of 260 DAEs). Note that max time step was chosen to be 0.05, otherwise the code might run even faster.

Example 5: This example was chosen to show the scalability of the solver for two-dimensional PDEs. The PDE system is simulated with a cell-centered finite difference method in x. The implemented solver is able to simulate this model in 37 time-steps and CPU time of 3.0 s with the IMPTRAP method with 0 failed time-steps with N×M= 64×64 elements (total of 4352 DAEs). Note that the maximum time step was chosen to be 0.25, otherwise the code might run even faster.

Example 6: This example was chosen to show the scalability of the solver for two-dimensional coupled PDEs resulting in large-scale DAEs and similarity of CPU time requirements with ODE models (example 5). The model chosen has a singularity at x = 0 and y = 0.5 and the model needs a large number of grids for uniform node spacing. The model solves both concentration and potential in an electrolyte. The PDE system is simulated with a cell-centered finite difference method in x. The implemented solver is able to simulate this model in 38 time-steps and CPU time of 16.65 s with the IMPTRAP method with 0 failed time-steps with N×M= 64×128 grids (total of 17152 DAEs). Note that the maximum time step was chosen to be 0.05 and the maximum growth rate was restricted to 3 (in the solver), otherwise the code might run even faster.

Table 2. CPU time for different algorithms for different examples considered. For all the models chosen, the following solver parameters were fixed $h_{init} = 10^{-6}$, $a_{tol}=10^{-6}$.

| Example | Methods | CPU Time [ms, s, m] | Number of DAEs | Number of Time Steps | Number of Rejected Steps |
|---|---|---|---|---|---|
| 1 | CN | 33,31,26 | 2 | 31 | 0 |
| 1 | IMPTRAP | 33,32,32 | 2 | 31 | 0 |
| 1 | EB | 81,84,79 | 2 | 129 | 1 |
| 1 | Rad | 31,30,31 | 2 | 31 | 0 |
| 2 | CN | 169,165,206 | 2 | 299 | 9 |
| 2 | IMPTRAP | 153,158,160 | 2 | 278 | 17 |
| 2 | EB | 1.26,1.32,1.27 | 2 | 2687 | 3 |
| 2 | Rad | 97,131,94 | 2 | 114 | 17 |
| 3 | CN | 49,47,49 | 2 | 55 | 0 |
| 3 | IMPTRAP | 39,47,56 | 2 | 55 | 0 |
| 3 | EB | 158,157,166 | 2 | 321 | 0 |
| 3 | Rad | 41,32,34 | 2 | 34 | 0 |
| 4 (N=128) | CN | 167,161,170 | 260 | 36 | 0 |
| | IMPTRAP | 112,161,165 | | 36 | 0 |
| | EB | 266,276,314 | | 74 | 4 |
| | Rad | 273,312,264 | | 40 | 4 |
| 4 (N=256) | CN | 259,253,289 | 516 | 35 | 0 |
| | IMPTRAP | 244,283,292 | | 35 | 0 |
| | EB | 480,488,484 | | 64 | 3 |
| | Rad | 626,642,645 | | 44 | 6 |
| 5 (N=64, M=64) Phi = 0.5 UMFPACK PARDISO | CN | 2.35,2.33,2.48 3.52,3.66,3.62 | 4352 | 37 | 0 |
| | IMPTRAP | 2.12,2.09,2.10 2.96,2.48,3.00 | | 37 | 0 |
| | EB | 2.44,2.43,2.47 3.61,3.72,3.55 | | 40 | 3 |
| | Rad | 5.07,5.10,5.11 7.03,6.91,7.08 | | 41 | 4 |
| 5 (N=128, M=128) Phi = 0.5 UMFPACK PARDISO | CN | 12.39,12.11,12.23 15.98,16.26,15.90 | 16896 | 40 | 2 |
| | IMPTRAP | 11.60,11.58,11.60 14.12,14.59,14.48 | | 40 | 2 |
| | EB | 12.73,12.68,12.48 16.82,16.61,16.57 | | 43 | 6 |
| | Rad | 25.78,25.95,25.70 32.26,32.16,31.95 | | 42 | 5 |
| 6 (N=128, M=256) UMFPACK PARDISO | CN | 83.21,84.50,84.98 79.99,85.22,90.16 | 67072 | 38 | 1 |
| | IMPTRAP | 83.59,84.11,83.02 88.93,82.38,83.37 | | 38 | 1 |
| | EB | 102.27,100.84,108.84 105.47,105.06,102.44 | | 43 | 6 |
| | Rad | 3.98,4.01,4.19 3.40,3.38,3.35 | | 43 | 4 |
| 6 (N=256, M=512) UMFPACK PARDISO | CN | 8.68,9.14,8.65 7.79,7.84,7.80 | 265216 | 36 | 0 |
| | IMPTRAP | 8.13,9.15,8.17 7.46,7.65,7.31 | | 36 | 0 |
| | EB | 10.35,10.75,10.29 9.21,8.95,8.97 | | 42 | 6 |
| | Rad | 41.33,43.56,41.37 23.48,22.73,22.50 | | 42 | 4 |

While it is not a fair comparison, MATLAB's ode15i was also used to benchmark the results obtained and the CPU time for different models. MATLAB's ode15i was run with maximum order of 2 to ensure A stability and with the same solver parameters for absolute tolerance, relative tolerance and maximum time step. This was done without providing the sparsity pattern or analytic

Jacobian forcing MATLAB to find the numerical Jacobian. MATLAB's ode15i was found to be similar in efficiency (often times more efficient) compared to the developed solver for the same solver parameters even though it calculates the numerical Jacobian for problems 1-3. When the problem size increases to more than 1000 DAEs (examples 4-6), the developed solver is more efficient. In particular for example 6, the MATLAB's ode15i cannot handle more than 40x80 grid points (it even required more than one hour of simulation time). Of course, optimizing the MATLAB code and providing the sparse analytic Jacobian can alleviate this. However, to use the developed DAE solver in Maple, the user has to provide only the DAEs, neither the sparsity pattern, nor the analytic Jacobian. MATLAB codes that solve examples 1-6 using ode15i can be obtained upon request from the corresponding author.

## 7. Future work, Perspectives, Code-dissemination, and Summary

In this paper, a sparse DAE solver was developed and implemented in Maple. Some of the aspects of the solver and future work are summarized below,

1. Four different algorithms were considered and implemented in this paper. For largescale problems from the discretization of PDEs in 2D, IMPTRAP method was found to be the most efficient. For extremely stiff problems, both Euler-backward and Radau methods might be more efficient. In particular, if unrealistic oscillations are observed, Euler-backward method is recommended. The code finds the analytic Jacobian by differentiating the functions in the model equations. If the function is not differentiable, then the code may not work. The solver has been found to work for some discontinuous cases. For example 1, replacing right-hand side in equation 1 with a piecewise function in z as given below works.

$$\frac{d}{dt} y(t) = z(t) \left( \begin{cases} 1 & 0.700000000000000 \leq z(t) \\ \frac{1}{2} & \text{otherwise} \end{cases} \right)$$

2. Continuous extension in time was not included in the current version of the solver. Maple's spline functions can be used for this. Our experience suggests that interpolations and continuous extensions are not as accurate and stable for the algebraic variables as the predicted values at the terminals (end-of-time steps) in the algorithms. Also, when a particular spatial discretization is used in a model with singularities (as in example 6), not all the variables will converge at the same order of accuracy requiring different orders for interpolation and continuous extension for different variables.
3. Implicit DAEs and ODEs (problems with non-constant mass-matrix form) can be addressed by making small changes in the algorithms, in particular for the Radau method. Other popular algorithms like BDF2, Radau5, TRBDF2, MEBDF will be implemented in the future.
4. Future versions of the solver include integration with Krylov type linear solvers and parallel sparse direct solvers,[9] event detection, *etc*. Reordering the equations and variables can reduce the bandwidth of the linear solver and reduce the computational time further.

5. Use of hardware floats and compiled codes will make the codes run faster with minimal memory requirements. This is easily doable based on Gauss elimination methods already implemented in Maple for 100-200 DAEs, however this will not scale for a large number of DAEs. Providing the residues and the Jacobian in an efficient code/format might speed up the code even further. This requires minor changes to the code to take the function F and Jacobian Jac as inputs to the solver.
6. The code is developed and distributed under MIT license without any restrictions. Two versions of the developed DAE solver are provided. The UMFPACK version (the default solver DAESolver.txt) works in any Maple installations by directly calling from Maple. The user should have a valid Intel MKL license and the required binaries to use the PARDISO-based DAE solver[11-13], DAESolverP.txt which is also directly called from Maple.

## Acknowledgments

Funds from Texas Materials Institute and UT's start-up package and past funds from the Department of Energy and multiple industries are gratefully acknowledged. The adaptive time-stepping procedure implemented in this code is taught by VS as a part of his graduate-level numerical methods course. Inputs from various students at UT and UW who took VS's graduate course on numerical methods are gratefully acknowledged.